\documentclass[11pt]{article}
\usepackage{amsmath,amsfonts,amssymb,latexsym,amsthm,dsfont}
\usepackage{geometry,graphicx}
\usepackage[backref]{hyperref}
\usepackage[latin1]{inputenc}
\usepackage[T1]{fontenc}

\newcommand{\ind}{\mathds{1}}
\newtheorem{thm}{Theorem}[section]%
\newtheorem{prop}[thm]{Proposition}%
\newtheorem{rem}[thm]{Remark}%
\newcommand{\dE}{\mathbb{E}}

\newcommand{\dN}{\mathbb{N}}
\newcommand{\dP}{\mathbb{P}}
\newcommand{\dR}{\mathbb{R}}

\newcommand{\cA}{\mathcal{A}}

\newcommand{\cL}{\mathcal{L}}
\newcommand{\cM}{\mathcal{M}}

\newcommand{\ABS}[1]{{{\left| #1 \right|}}} 
\newcommand{\BRA}[1]{{{\left\{#1\right\}}}} 
\newcommand{\PAR}[1]{{{\left(#1\right)}}} 
\renewcommand{\leq}{\leqslant}
\renewcommand{\geq}{\geqslant}

\title{On the Laplace transform of perpetuities with thin tails}

 \author{%
  Jean-Baptiste~\textsc{Bardet}, %
  H\'el\`ene~\textsc{Gu\'erin}, %
  Florent~\textsc{Malrieu}}

\date{Unpublished note -- \today}

\begin{document}

\maketitle

\begin{abstract}
We consider the random variables $R$ which are solutions of the distributional 
equation $R\overset{\cL}{=}MR+Q$, where $(Q,M)$ is independent of $R$ and 
$\ABS{M}\leq 1$. Goldie and Gr\"ubel showed that the tails of $R$ are no heavier 
than exponential. Alsmeyer and \emph{al} provide a complete description of the domain 
of the Laplace transform of $R$. We present here a simple proof in a particular case 
and an extension to the Markovian case. 
\end{abstract}

\noindent\textbf{AMS Classification 2000: }Primary 60H25; secondary 60E99 

\section{Introduction}

We define on some probability space $(\Omega, \cA,\dP)$ a couple of random variables $(M,Q)$, 
a sequence ${(M_n,Q_n)}_{n\geq 0}$ of independent and identically distributed  random vectors
with the same law as $(M,Q)$, and $R_0$ a random variable independent of the sequence ${(M_n,Q_n)}_{n\geq 0}$.
Define the sequence ${(R_n)}_{n\geq 0}$ by
\begin{equation}\label{eq:rn}
R_{n+1}=M_nR_n+Q_n,
\end{equation}
for any $n\geq 0$. This sequence has been extensively studied in the last decades. Under weak assumptions (see \cite{vervaat}) which are obviously fullfilled in our setting, it can be shown that the sequence ${(R_n)}_{n\geq 0}$ converges almost surely to a random variable $R$ such that 
\begin{equation}\label{eq:autoR}
R \overset{\cL}{=}M R+Q,
\end{equation}
where $R$ is independent of $(M,Q)$.

In \cite{kesten}, Kesten established that $R$ is in general heavy-tailed (\emph{i.e.} not all the moments of $R$ are finite) even if $Q$ is light-tailed as soon as $\ABS{M}$ can be greater than 1. Nevertheless, Goldie and Gr\"ubel \cite{goldie-grubel} have shown that $R$ can have some exponential moments if $\ABS{M}\leq 1$. In particular, if $Q$ and $M$ are nonnegative the following result holds. 

\begin{thm}[Goldie, Gr\"ubel \cite{goldie-grubel}]\label{th:GG}
Assume that 
$$
\dP(Q\geq0, 0\leq M\leq 1)=1,\quad \dP(M<1)>0
$$
 and that there is $v_Q>0$ 
(possibly infinite) such that 
\begin{equation}\label{eq:vQ}
\dE\PAR{e^{vQ}}
\begin{cases}
 < +\infty & \text{if $v<v_Q$,}\\
 = +\infty & \text{if $v>v_Q$.}
\end{cases}
\end{equation}
Then, the Laplace transform $v\mapsto \dE\PAR{e^{vR}}$ of the solution $R$ of \eqref{eq:autoR} is finite on the set 
$(-\infty,v_{GG})$ with $v_{GG}=v_Q\wedge\sup\BRA{v\geq 0,\ \dE(e^{vQ}M)<1}$.
\end{thm}
In fact, the domain of the Laplace transform of $R$ is larger than $(-\infty,v_{GG})$. In \cite{iksanov}, a full 
description of this domain is established. Let us provide a simple proof under the assumptions of Theorem \ref{th:GG}.

\section{The main result}

 \begin{thm}\label{th:GGoptimal}
Under the assumptions of Theorem \ref{th:GG}, assuming furthermore 
that $R_{0}$ is non-negative and has all
its exponential moments finite, then
$$
\sup_{n\geq0}\dE\PAR{e^{vR_n}}<+\infty%
\quad\text{and}\quad
\dE\PAR{e^{vR}}<+\infty
$$
for any $v<v_c$, where 
$$
v_c=v_Q\wedge \sup\{v\geq0,\,\dE\PAR{e^{vQ}1_{\{M=1\}}}<1\}.
$$

Moreover, for any  $v>v_c$, $\sup_{n\geq0}\dE\PAR{e^{vR_n}}=+\infty$ and $\dE\PAR{e^{vR}}=+\infty$.
\end{thm}

For other recent generalizations of \cite{goldie-grubel}, the interested reader is 
referred to \cite{hitwes}, where the authors give sharper results than ours on the 
tails for some specific examples.

\begin{proof}[Proof of Theorem \ref{th:GGoptimal}]
 Let us start this section with the main lines of the 
proof of Theorem \ref{th:GG} of Goldie and Gr\"ubel \cite{goldie-grubel}. For $\rho>0$, 
let $\cM_\rho$ be the set of probability measures on $\dR_{+}$ 
 with finite exponential moment of order $\rho$, and $d_\rho$ a distance defined on $\cM_\rho$ by: 
for $\mu,\nu\in\cM_\rho$, 
 $$
 d_\rho(\mu,\nu)=\int_0^\infty\!e^{\rho u}\ABS{\mu[u,\infty)-\nu[u,\infty)}\,du.
 $$
Define the application $T$ on $\cM_\rho$ as follows: for $X$ with law $\mu\in \cM_\rho$, $T\mu$ is the law of $Q+MX$ with $(M,Q)$ independent of $X$. It is shown in \cite{goldie-grubel} that,  
  $$
 d_\rho(T\mu,T\nu)\leq \dE(e^{\rho Q}M) d_\rho(\mu,\nu). 
 $$
  Since 
 $$
 \dE\PAR{e^{vX}}=v\int_0^\infty\!e^{vu}\dP(X\geq u)\,du,
 $$
 one can show that, for any $n\geq 0$ and $v<\min(v_0,v_Q)$ with $v_0=\sup\{v\geq0,\,\dE\PAR{e^{vQ}M}<1\}$, 
\begin{align*}
\dE\PAR{e^{vR_n}}&\leq v \frac{1-\dE(e^{vQ}M)^n}{1-\dE(e^{vQ}M)}d_v(T\mu_0,\mu_0)+\dE\PAR{e^{vR_0}}.
\end{align*}
In others words, Goldie and Gr\"ubel \cite{goldie-grubel} established that for any $v<\min(v_0,v_Q)$,
${(\dE(e^{vR_n}))}_{n}$ is uniformly bounded. This estimate can be extended to a larger domain. 

Let us define $v_1=\sup\{v\geq0,\,\dE\PAR{e^{vQ}1_{\{M=1\}}}<1\}$ and $v_c=\min(v_1,v_Q)$.
Let us fix $v<v_c$ and choose $\varepsilon>0$ such that 
 $$
 \rho:=\dE(e^{vQ}1_{\{1-\varepsilon<M\leq1\}})<1.
 $$
  Then we get, for any $n\geq0$,
\begin{align*}
 L_{n+1}(v):=\dE(e^{vR_{n+1}})&=\dE(e^{v(M_nR_n+Q_n)})\\
 &\leq \dE(e^{v((1-\varepsilon)R_n+Q_n)}1_{\{M_n\leq 1-\varepsilon\}})+\dE(e^{v(R_n+Q_n)}1_{\{1-\varepsilon<M_n\leq 1\}})\\
 &\leq L_n((1-\varepsilon)v)L_Q(v)+\rho L_n(v)
\end{align*}
where $L_Q(v)=\dE(e^{vQ})$. By iteration of this estimate, one gets for any $n\geq0$
$$
L_n(v)\leq \Big(\sum_{k=0}^{n-1}\rho^kL_{n-k}((1-\varepsilon)v)\Big)L_Q(v)+\rho^nL_0(v)\,.
$$
Let us notice that we have in fact more: for the same $\varepsilon$, and for any $\tilde{v}\leq v$,
$\tilde{\rho}:=\dE(e^{\tilde{v}Q}1_{\{1-\varepsilon<M\leq1\}})<\rho$, hence, by the same method 
as before,
\begin{equation}
\label{eq:iteratelaplace}
L_n(\tilde{v})\leq \Big(\sum_{k=0}^{n-1}\rho^kL_{n-k}((1-\varepsilon)\tilde{v})\Big)L_Q(\tilde{v})+\rho^nL_0(\tilde{v})\,.
\end{equation}

Let us define $\overline{L}=\sup_{n\geq0} L_n$. Taking the supremum over $n$ in \eqref{eq:iteratelaplace}, one gets for any $\tilde{v}\leq v$
\begin{equation}
\label{eq:iteratelaplacesup}
\overline{L}(\tilde{v})\leq \frac{1}{1-\rho}\overline{L}((1-\varepsilon)\tilde{v})L_Q(\tilde{v})+L_0(\tilde{v})\,.
\end{equation}
There is $k\in\dN$ such that ${(1-\varepsilon)}^kv<v_0$, hence $ \overline{L}({(1-\varepsilon)}^kv)<+\infty$. Applying $k$ times estimate \eqref{eq:iteratelaplacesup}, one then obtains immediatly that $ \overline{L}(v)<+\infty$, which achieves the first part of the proof.

On the other hand, if $v>v_Q$, $R_1\geq Q_0$ immediatly implies $L_1(v)=+\infty$; if $v>v_1$, 
$\rho_0:=\dE(e^{vQ}1_{\{M=1\}})>1$ (except in a trivial case, left to the reader) and, for all $n\geq0$,
$$
L_{n+1}(v)\geq \rho_0 L_n(v)\,,$$
implying that $\overline{L}(v)=+\infty$.
\end{proof}

\section{Some extensions and perspectives}

What happens if the random variables ${(M_n,Q_n)}_{n\geq 0}$ are no longer independent? We provide here a partial result under a Markovian assumption when the contractive term $M$ is less than 1.

Let us introduce $X=(X_n)_{n\geq 0}$ an irreducible recurrent Markov process with finite space $E$ and $\PAR{(M_n(x),Q_n(x))_{x\in E}}_{n\geq 0}$ a sequence of i.i.d. random vectors supposed to be independent of $X$. We assume that, for all $x\in E$,
$$
\dP\PAR{0\leq M(x)< 1}=1,
$$
but we do not assume in the sequel that $Q$ is non negative. The sequence $(R_n)_{n\geq 0}$ is defined by 
$$
 R_{n+1}=M_n(X_n)R_n+Q_n(X_n),
 $$
$R_0$ being arbitrary (with all exponential moments).
Notice that the process $(X_n,R_n)_{n\geq0}$ is a Markov process whereas $(R_{n})_{n\geq0}$ 
is not (in general). 

\begin{prop}\label{prop:markovmajoration}
Introduce $\underline v=\inf_{x\in E}v_{\ABS{Q(x)}}$, with $v_{\ABS{Q(x)}}$ defined as in (\ref{eq:vQ}). For any $v< \underline v$, 
 $$
 \sup_{n\geq 0}\dE\PAR{e^{v\ABS{R_n}}}<+\infty. 
 $$
 Moreover, if $v>\underline v$, then this supremum is infinite.
\end{prop}

\begin{proof}
Let us introduce $\overline M_n=\max_{x\in E}M_n(x)$ and $\overline Q_n=\max_{x\in E}\ABS{Q_n(x)}$. The random variables ${((\overline M_n,\overline Q_n))}_{n\geq 0}$ are i.i.d. Define the sequence ${(\overline R_n)}_{n\geq 0}$ by
$$
\overline R_0=\ABS{R_0}%
\quad\text{and} \quad%
 \overline{R}_{n+1}=\overline M_n\overline{R}_n+\overline Q_n
 \quad\text{for } n\geq 1.
 $$
Obviously, $\ABS{R_n}\leq \overline{R}_n$ for all $n\geq 0$. Thus it is sufficient to study the Laplace transforms of $({\overline R_n})_{n\geq 0}$. On the other hand, Theorem \ref{th:GGoptimal} ensures that  
${(\dE\PAR{e^{v\overline R_n}})}_n$ is uniformly bounded
as soon as $v<\overline v_c=\min(\overline v_1,v_{\overline Q})$ with $\overline v_1=\sup\{v\geq0\,:\,\dE(e^{v\overline Q}1_{\{\overline M=1\}})<1\}>0$. In our case, $\overline v_1$ is infinite since $\dP(\overline M<1)=1$. At last, for $v\geq0$,
$$
\sup_{x\in E}\dE\PAR{e^{v\ABS{Q(x)}}}\leq%
\dE\PAR{e^{v\overline Q}}=\dE\PAR{\sup_{x\in E}e^{v\ABS{Q(x)}}}%
\leq \sum_{x\in E}\dE\PAR{e^{v\ABS{Q(x)}}}.
$$
Thus $v_{\overline Q}=\inf_{x\in E}v_{\ABS{Q(x)}}$. 

On the other hand, choose $v>\underline v$. There exists $x_0\in E$ such that $\dE\PAR{e^{v\ABS{Q(x_0)}}}$ is infinite. Then, for any $n\geq 0$,
\begin{align*}
\dE\PAR{e^{v\ABS{R_{n+1}}}}&\geq \dE\PAR{e^{v\ABS{R_{n+1}}}\ind_\BRA{X_n=x_0}}\\
&\geq \dE\PAR{e^{-v\ABS{R_n}}e^{v\ABS{Q_n(x_0)}}\ind_\BRA{X_n=x_0}}\\
&\geq \dE\PAR{\ind_\BRA{X_n=x_0}e^{-v\ABS{R_n}}}\dE\PAR{e^{v\ABS{Q_n(x_0)}}}.
 \end{align*}
The recurrence of $X$ ensures that $\BRA{n\geq 0,\ \dE\PAR{e^{v\ABS{R_n}}}=+\infty}$ is infinite. 
\end{proof}

\begin{rem}
In \cite{nous}, we use the previous estimates to improve the results of \cite{guyon,saporta-yao} on the tails of the invariant measure of a diffusion process with Markov switching. 
\end{rem}

\addcontentsline{toc}{section}{\refname}%
{ \footnotesize
\bibliography{notebgm}
}
\bibliographystyle{amsplain}

\bigskip

\begin{flushright}\texttt{Compiled \today.}\end{flushright}


{\footnotesize %
 \noindent Jean-Baptiste \textsc{Bardet}
e-mail: \texttt{jean-baptiste.bardet(AT)univ-rouen.fr}
  
 \medskip

 \noindent\textsc{UMR 6085 CNRS Laboratoire de Math\'ematiques Rapha\"el Salem (LMRS)\\
Universit\'e de Rouen,
Avenue de l'Universit\'e, BP 12,
F-76801 Saint Etienne du Rouvray}

 \bigskip
  
 \noindent H\'el\`ene \textsc{Gu\'erin},
e-mail: \texttt{helene.guerin(AT)univ-rennes1.fr}

 \medskip

  \noindent\textsc{UMR 6625 CNRS Institut de Recherche Math\'ematique de
    Rennes (IRMAR) \\ Universit\'e de Rennes I, Campus de Beaulieu, F-35042
    Rennes \textsc{Cedex}, France.}

 \bigskip
  
 \noindent Florent \textsc{Malrieu}, corresponding author,
 e-mail: \texttt{florent.malrieu(AT)univ-rennes1.fr}

 \medskip

  \noindent\textsc{UMR 6625 CNRS Institut de Recherche Math\'ematique de
    Rennes (IRMAR) \\ Universit\'e de Rennes I, Campus de Beaulieu, F-35042
    Rennes \textsc{Cedex}, France.}

}

\end{document}